\documentclass{article}
\usepackage{graphicx}
\usepackage{amssymb}
\usepackage{amsmath}

\pdfoutput=1

\input{tcilatex}
\begin{document}

\title{ACCELERATION\ CONTROL\ IN\ NONLINEAR\ VIBRATING\ SYSTEMS\ BASED\ ON\
DAMPED\ LEAST\ SQUARES\ }
\author{V.N. Pilipchuk \\
Wayne State University\\
Detroit, MI 48202}
\maketitle

\begin{abstract}
A discrete time control algorithm using the damped least squares is
introduced for acceleration and energy exchange controls in nonlinear
vibrating systems. It is shown that the damping constant of least squares
and sampling time step of the controller must be inversely related to insure
that vanishing the time step has little effect on the results. The algorithm
is illustrated on two linearly coupled Duffing oscillators near the 1:1
internal resonance. In particular, it is shown that varying the dissipation
ratio of one of the two oscillators can significantly suppress the nonlinear
beat phenomenon.
\end{abstract}

\section{Introduction}

The damped least squares is a simple but effective analytical manipulation
that helps to avoid singularity in practical minimization and control
algorithms. It is also known as Levenberg-Marquardt method \cite%
{Levenberg:1944}. In order to illustrate the idea in simple terms, let us
consider the minimization problem 
\begin{equation}
\left\Vert E-A\delta u\right\Vert ^{2}\rightarrow \min  \label{P0}
\end{equation}%
where $E$ $\in R^{n}$ is a given vector, the notation $\left\Vert
...\right\Vert $ indicates the Euclidean norm in $R^{n}$, $A$ is typically a
Jacobian matrix of $n$ rows and $m$ columns, and $\delta u\in R^{m}$ is an
unknown minimization vector. Although a formal solution of this problem is
given by $\delta u=(A^{T}A)^{-1}A^{T}E$, the matrix product $A^{T}A$ may
appear to be singular so that no unique solution is possible. This fact
usually points to multiple possibilities of achieving the same result unless
specific conditions are imposed on the vector $\delta u$. The idea of damped
least squares is to avoid such conditioning by adding one more quadratic
form to the left hand side of expression (\ref{P0}) as follows%
\begin{equation}
\left\Vert E-A\delta u\right\Vert ^{2}+\lambda \left\Vert \delta
u\right\Vert ^{2}\rightarrow \min  \label{P}
\end{equation}%
where $\lambda $ is a positive scalar number, which is often called \textit{%
damping constant}; note that the term `damping' has no relation to the
physical damping or energy dissipation effects in vibrating systems usually
characterized by \textit{damping ratios}.

Now the inverse matrix includes the damping constant $\lambda $ which can
provide the uniqueness of solution given by

\begin{equation}
\delta u=(A^{T}A+\lambda I)^{-1}A^{T}E  \label{du}
\end{equation}%
where $I$ is $n\times n$ identity matrix.

Different arguments are discussed in the literature regarding the use of
damped least squares and best choice for the damping parameter $\lambda $ 
\cite{Chan:1988}, \cite{Chiaverini:1988}, \cite{Chiaverini:1994}, \cite%
{Chung:2000}, \cite{Deo:1992}, \cite{Deo:1993}, \cite{Hamada:2006}, \cite%
{Hattori:2003}, \cite{Mayorga:1990}, \cite{Mayorga:1992}, \cite%
{Nakamura:1986JDSMC}, \cite{Wampler:1986IEEE}, \cite{Wampler:1988JDSMC}. In
particular, it was noticed that the parameter $\lambda $ may affect
convergence properties of the corresponding algorithms. The parameter $%
\lambda $ can be used also for other reason such as shifting the solution $%
\delta u$ into desired area in $R^{m}$. \ In this case, the meaning of $%
\lambda $ is rather close to that of Lagrangian multiplier imposing
constraints on control inputs.

In case of dynamical systems, when all the quantities in (\ref{P}) may
depend on time, a continuous time analogue of (\ref{P}) can be written in
the integral form%
\begin{equation}
\min_{\delta u}\int_{0}^{T}(\left\Vert E-A\delta u\right\Vert ^{2}+\lambda
\left\Vert \delta u\right\Vert ^{2})dt  \label{P_cont}
\end{equation}%
where the interval of integration is manipulated as needed, for instance, $T$
can be equal to sampling time of the controller \cite{Lee:2009}.

However, in the present work, a discrete time algorithm based on the damped
least squares solution (\ref{du}), which is used locally at every sample
time $t_{n}$, is introduced. \ Such algorithm appears to be essentially
discrete namely using different time step $h$ may lead to different results.
Nevertheless, if the parameters $\lambda $ and $h$ are coupled by some
condition then the control input and system response show no significant
dependence on the time step.

A motivation for the present work is as follows. In order to comply with the
standard tool of dynamical systems dealing with differential equations, the
methods of control are often formulated in continuous time by silently
assuming that a discrete time analogous is easy to obtain one way or another
whenever it is needed for practical reasons. For instance, data acquisition
cards and on-board computers of ground vehicles usually acquire and process
data once per 0.01 sec. Typically, based on the information, which is known
about the system dynamic states and control inputs by the time instance $%
t_{n}$, the computer must calculate control adjustments for the next active
time instance, $t_{n+1}$. The corresponding computational time should not
therefore exceed $t_{n+1}-t_{n}=0.01$ sec. Generally speaking, it is
possible to memorize snapshots of the dynamic states and control inputs at
some of the previous times $\{$..., $t_{n-2}$, $t_{n-1}\}$. However,
increasing the volume of input data may complicate the code and, as a
result, slow down the calculation process. Therefore, let us assume that
updates for the control inputs are obtained by processing the system states,
controls, and target states given only at the current time instance, $t_{n}$%
. The corresponding algorithm can be built on the system model described by
its differential equations of motion and some rule for minimizing the
deviation (error) of the current dynamic states from the target. Recall
that, in the present work, such a rule will be defined according to the
damped least squares (\ref{P}). \ Illustrating physical example of two
linearly coupled Duffing oscillators is considered. It is shown that the
corresponding algorithm, which is naturally designed and effectively working
in discrete time, may face a problem of transition to the continuous time
limit.

\section{Problem formulation}

\noindent Consider the dynamical system 
\begin{equation}
\ddot{x}=f(x,\dot{x},t,u)  \label{system}
\end{equation}%
where $x=x(t)\in R^{n}$ is the system position (configuration) vector, the
overdot indicates derivative with respect to time $t$, the right-hand side $%
f\in R^{n}$ represents a vector-function that may be interpreted as a force
per unit mass of the system, and $u=u(t)\in R^{m}$ is a control vector,
whose dimension may differ from that of the positional vector $x$ so that
generally $n\neq m$.

In common words, the purpose of control $u(t)$ is to keep the acceleration $%
\ddot{x}(t)$ of system (\ref{system}) as close as possible to the target $%
\ddot{x}^{\ast }(t)$. The term `close' will be interpreted below through a
specifically designed target function of the following error vector 
\begin{equation}
E(t)=\ddot{x}^{\ast }(t)-\ddot{x}(t)  \label{error cont}
\end{equation}

As discussed in Introduction, for practical implementations, the problem
must be formulated in terms of the discrete time $\{t_{k}\}$ as follows. Let 
$x_{k}=x(t_{k})$, $\dot{x}_{k}=\dot{x}(t_{k})$, and $u_{k}=u(t_{k})$ are
observed at some time instance $t_{k}$. The corresponding target
acceleration, $\ddot{x}_{k}^{\ast }=\ddot{x}^{\ast }(t_{k})$, is assumed to
be known. Then, taking into account (\ref{system}) and (\ref{error cont}),
gives the following error at the same time instance%
\begin{equation}
E_{k}=\ddot{x}_{k}^{\ast }-f(x_{k},\dot{x}_{k},t_{k},u_{k})
\label{error discr}
\end{equation}

Now the purpose of control is to minimize the following target function%
\begin{eqnarray}
P_{k} &=&\frac{1}{2}E_{k}^{T}W_{k}E_{k}  \label{target 1} \\
&=&\frac{1}{2}[\ddot{x}_{k}^{\ast }-f(x_{k},\dot{x}%
_{k},t_{k},u_{k})]^{T}W_{k}[\ddot{x}_{k}^{\ast }-f(x_{k},\dot{x}%
_{k},t_{k},u_{k})]  \notag
\end{eqnarray}%
where $W_{k}$ is $n\times n$ diagonal weight matrix whose elements are
positive or at least non-negative functions of the system states, $%
W_{k}=W(x_{k},\dot{x}_{k},t_{k})$.

Note that all the quantities in expression (\ref{target 1}) represent a
snapshot of the system at $t=t_{k}$ while including no data from the
previous time step $t_{k-1}$. Since the control vector $u_{k}$ cannot be
already changed at time $t_{k}$ then quantity $P_{k}$ is out of control at
time $t_{k}$. In other words expression (\ref{target 1}) summarizes all what
is observed \textit{now}, at the time instance $t_{k}$. The question is how
to adjust the control vector $u$ for the next step $t_{k+1}$ based on the
information included in (\ref{target 1}) while the system state at $%
t=t_{k+1} $ is yet unknown, and no information from the previous times $\{$%
..., $t_{n-2}$, $t_{n-1}\}$ is available.

Let us represent such an update for the control vector in the form%
\begin{equation}
u_{k+1}=u_{k}+\delta u_{k}  \label{u(n+1)}
\end{equation}%
were $\delta u_{k}$ is an unknown adjustment of the control input.

Replacing $u_{k}$ in (\ref{target 1}) by (\ref{u(n+1)}) and taking into
account that%
\begin{eqnarray}
f(x_{k},\dot{x}_{k},t_{k},u_{k+1}) &=&f(x_{k},\dot{x}_{k},t_{k},u_{k})+A_{k}%
\delta u_{k}+O(\left\Vert \delta u_{k}\right\Vert ^{2})  \label{upredict} \\
A_{k} &=&\partial f(x_{k},\dot{x}_{k},t_{k},u_{k})/\partial u_{k}  \notag
\end{eqnarray}%
gives%
\begin{equation}
P_{k}=\frac{1}{2}(E_{k}-A_{k}\delta u_{k})^{T}W_{k}(E_{k}-A_{k}\delta u_{k})
\label{target 2}
\end{equation}%
where $A_{k}$ is the Jacobian matrix of $n$ rows and $m$ columns.

Although the replacement $u_{k}$ by $u_{k+1}$ in (\ref{upredict}) may look
artificial, this is how the update rule for the control vector $u$ is
actually defined here. Namely, if $u_{k}$ did not provide a minimum for $%
P_{k}(\ddot{x}_{k}^{\ast },x_{k},\dot{x}_{k},t_{k},u_{k})$, then let us
minimize $P_{k}(\ddot{x}_{k}^{\ast },x_{k},\dot{x}_{k},t_{k},u_{k}+\delta
u_{k})$ with respect to $\delta u_{k}$ and then apply the adjusted vector (%
\ref{u(n+1)}) at least the next next time, $t_{n+1}$. Assuming that the
variation $\delta u_{k}$ is small, in other words, $u_{k}$ is still close
enough to the minimum, expansion (\ref{upredict}) is applied. Now the
problem is formulated as a minimization of the quadratic form (\ref{target 2}%
) with respect to the adjustment $\delta u_{k}$. However, what often happens
practically is that function (\ref{target 2}) has no unique minimum so that
equation%
\begin{equation}
\frac{dP_{k}}{d\delta u_{k}}=0  \label{equation du}
\end{equation}%
has no unique solution. In addition, even if the unique solution does exist,
it may not satisfy some conditions imposed on the control input due to the
physical specifics of actuators. As a result, some constraint conditions may
appear to be necessary to impose on the variation of control adjustment, $%
\delta u_{k}$. However, the presence of constraints would drastically
complicate the problem. Instead, the target function (\ref{target 2}) can be
modified in order to move solution $\delta u_{k}$ into the allowed domain.
For that reason, let us generalize function (\ref{target 2}) as\ \ 
\begin{eqnarray}
P_{k} &=&\frac{1}{2}(E_{k}-A_{k}\delta u_{k})^{T}W_{k}(E_{k}-A_{k}\delta
u_{k})  \notag \\
&&+\frac{1}{2}(B_{k}+C_{k}\delta u_{k})^{T}\Lambda _{k}(B_{k}+C_{k}\delta
u_{k})  \label{target 3}
\end{eqnarray}%
where $\Lambda _{k}=\Lambda (x_{k},\dot{x}_{k},t_{k})$ is a diagonal
regularization matrix, $B_{k}=B(x_{k},\dot{x}_{k},t_{k})$ is a
vector-function of $n$ elements, and $C_{k}=C(x_{k},\dot{x}_{k},t_{k})$ is a
matrix of $n$ rows and $m$ columns.

Note that the structure of new function (\ref{target 3}) is a generalization
of (\ref{P}). Substituting (\ref{target 3}) in (\ref{equation du}), gives a
linear set of equations in the matrix form whose solution $\delta u_{k}$
brings relationship (\ref{u(n+1)}) to the form

\begin{equation}
u_{k+1}=u_{k}+(A_{k}^{T}W_{k}A_{k}+C_{k}^{T}\Lambda
_{k}C_{k})^{-1}(A_{k}^{T}W_{k}E_{k}-C_{k}^{T}\Lambda _{k}B_{k})
\label{soln u}
\end{equation}

The entire discrete time system is obtained by adding a discrete version of
the dynamical system (\ref{system}) to (\ref{soln u}) . Assuming that the
time step is fixed, $t_{k+1}-t_{k}=h$, a simple discrete version can be
obtained by means of Euler explicit scheme as follows%
\begin{eqnarray}
x_{k+1} &=&x_{k}+hv_{k}  \notag \\
v_{k+1} &=&v_{k}+hf(x_{k},v_{k},t_{k},u_{k})  \label{discrete system}
\end{eqnarray}

Finally, equations (\ref{soln u}) and (\ref{discrete system}) represent a
discrete time dynamical system, whose motion should follow the target
acceleration $\ddot{x}_{k}^{\ast }=\ddot{x}^{\ast }(t_{k})$.

It will be shown in the next section that the structure of equation (\ref%
{soln u}) does not allow for the transition to continuous limit of the
entire dynamic system (\ref{soln u}) through (\ref{discrete system}), unless
some specific assumption are imposed on the parameters in order to guarantee
that $\delta u_{k}=O(h)$ as $h\rightarrow 0$.

\section{The illustrating example}

The algorithm, which is designed in the previous section, is applied now to
a two-degrees-of-freedom nonlinear vibrating system for an \textit{active}
control of the energy exchange (nonlinear beats) between the two
oscillators. The problem of \ \textit{passive} control of energy flows in
vibrating systems is of great interest \cite{Vakakis:2009}, and it is
actively discussed from the standpoint of nonlinear beat phenomena \cite%
{ManevitchGendelman:2011}. The beating phenomenon takes place when
frequencies of the corresponding linear oscillators are either equal or at
least close enough to each other.

For illustrating purposes, let us consider two unit-mass Duffing oscillators
of the same linear stiffness $K$ coupled by the linear spring of stiffness $%
\gamma $. The system position is described by the vector-function of
coordinates, $x(t)=(x_{1}(t),x_{2}(t))^{T}$. Introducing the parameters $%
\Omega =(\gamma +K)^{1/2}$ \ and $\varepsilon =\gamma /(\gamma +K)$, brings
the differential equations of motion to the form 
\begin{eqnarray}
\dot{x}_{1} &=&v_{1}  \notag \\
\dot{x}_{2} &=&v_{2}  \notag \\
\dot{v}_{1} &=&-2\zeta \Omega v_{1}-\Omega ^{2}x_{1}+\varepsilon (\Omega
^{2}x_{2}-\alpha x_{1}^{3})\equiv f_{1}(x_{1},x_{2},v_{1})  \label{nb3} \\
\dot{v}_{2} &=&-2u\Omega v_{2}-\Omega ^{2}x_{2}+\varepsilon (\Omega
^{2}x_{1}-\alpha x_{2}^{3})\equiv f_{2}(x_{1},x_{2},v_{2},u)  \notag
\end{eqnarray}%
where $\alpha $ is a positive parameter, $\zeta $ and $u$ are damping ratios%
\footnote{%
As mentioned in Introduction, the damping (dissipation) ratio should not be
confused with the damping coefficient $\lambda $.} of the first and the
second oscillators, respectively; the damping ratio $u$, which is explicitly
shown as an argument of the function $f_{2}(x_{1},x_{2},v_{2},u)$, will be
considered as a control input.

The problem now is to find such variable damping ratio $u=u(t)$ under which
the second oscillator accelerates as close as possible to the given (target)
acceleration, $\ddot{x}_{2}^{\ast }(t)$.

Following the discussion of the previous section, let us consider the
problem in the discrete time $\{t_{k}\}$. In order to avoid confusion, the
iterator $k$ will be separated from the vector component indexes by coma,
for instance, $x_{k}=(x_{1,k},x_{2,k})^{T}$. Since only the second mass
acceleration is of interest and the system under consideration includes only
one control input $u$, then, assuming the weights to be constant, gives%
\begin{equation*}
W_{k}=\left[ 
\begin{array}{cc}
0 & 0 \\ 
0 & 1%
\end{array}%
\right] \text{, \ \ }A_{k}=\frac{\partial }{\partial u_{k}}\left[ 
\begin{array}{c}
f_{1,k} \\ 
f_{2,k}%
\end{array}%
\right]
\end{equation*}%
where $f_{1,k}\equiv f_{1}(x_{1,k},x_{2,k},v_{1,k})$ and $f_{2,k}\equiv
f_{2}(x_{1,k},x_{2,k},v_{2,k},u_{k})$, and other matrix terms become scalar
quantities, say, $\Lambda _{k}=\lambda $, \ \ $B_{k}=b$, \ \ and $C_{k}=1$.
The unities in $W_{k}$ and $C_{k}$ can always be achieved by re-scaling the
target function and parameters $\lambda $ and $b$. Note that re-scaling the
target function by a constant factor has no effect on the solution of
equation (\ref{equation du}).

As a result, the target function (\ref{target 3}) takes the form

\begin{equation}
P_{k}=\frac{1}{2}\left( \ddot{x}_{2,k}^{\ast }-f_{2,k}-\frac{\partial f_{2,k}%
}{\partial u_{k}}\delta u_{k}\right) ^{2}+\frac{\lambda }{2}(b+\delta
u_{k})^{2}  \label{sample_target}
\end{equation}

In this case, equation (\ref{equation du}) represents a single linear
equation with respect to the scalar control adjustment, $\delta u_{k}$.
Substituting the corresponding solution in (\ref{soln u}) and taking into
account (\ref{discrete system}), gives the discrete time dynamical system 
\begin{equation}
u_{k+1}=u_{k}-\frac{(f_{2,k}-\ddot{x}_{2,k}^{\ast })(\partial
f_{2,k}/\partial u_{k})+\lambda b}{(\partial f_{2,k}/\partial
u_{k})^{2}+\lambda }  \label{du_solution}
\end{equation}%
and%
\begin{eqnarray}
x_{1,k+1} &=&x_{1,k}+hv_{1,k}  \notag \\
x_{2,k+1} &=&x_{2,k}+hv_{2,k}  \notag \\
v_{1,k+1} &=&v_{1,k}+hf_{1,k}  \label{sample_iter} \\
v_{2,k+1} &=&v_{2,k}+hf_{2,k}  \notag
\end{eqnarray}

Let us assume now that the target acceleration $\ddot{x}_{2}^{\ast }$ is
zero, in other words, the purpose of control is to minimize acceleration of
the second oscillator at any sample time $t_{k}$ as much as possible. Let us
set still arbitrary parameter $b$ also to zero. Then the target function (%
\ref{sample_target}) and dynamical system (\ref{du_solution}) and (\ref%
{sample_iter}) take the form%
\begin{equation}
P_{k}=\frac{1}{2}\left[ f_{2}(x_{1,k},x_{2,k},v_{2,k},u_{k})+\frac{\partial
f_{2}(x_{1,k},x_{2,k},v_{2,k},u_{k})}{\partial u_{k}}\delta u_{k}\right]
^{2}+\frac{\lambda }{2}(\delta u_{k})^{2}  \label{sample_target1}
\end{equation}

\begin{eqnarray}
u_{k+1} &=&u_{k}+\frac{2\Omega v_{2,k}}{4\Omega ^{2}v_{2,k}^{2}+\lambda }%
f_{2}(x_{1,k},x_{2,k},v_{2,k},u_{k})  \notag \\
x_{1,k+1} &=&x_{1,k}+hv_{1,k}  \notag \\
x_{2,k+1} &=&x_{2,k}+hv_{2,k}  \label{sample_iter1} \\
v_{1,k+1} &=&v_{1,k}+hf_{1}(x_{1,k},x_{2,k},v_{1,k})  \notag \\
v_{2,k+1} &=&v_{2,k}+hf_{2}(x_{1,k},x_{2,k},v_{2,k},u_{k})  \notag
\end{eqnarray}%
where the functions $f_{1}$ and $f_{2}$ are defined in (\ref{nb3}).

As follows from the first equation in (\ref{sample_iter1}), transition to
the continuous time limit for the entire system (\ref{sample_iter1}) would
be possible under the condition that%
\begin{equation}
\frac{2\Omega v_{2,k}}{4\Omega ^{2}v_{2,k}^{2}+\lambda }=O(h)\text{, \ \ as
\ \ }h\rightarrow 0  \label{cond1}
\end{equation}

Condition (\ref{cond1}) can be satisfied by assuming that $\Omega =O(h)$.
Such an assumption, however, makes little if any physical sense. As an
alternative choice, the condition $\lambda =O(h^{-1})$ can be imposed by
setting, for instance, 
\begin{equation}
\lambda h=\lambda _{0}  \label{cond2}
\end{equation}%
where $\lambda _{0}$ remains finite as $h\rightarrow 0$.

However, condition (\ref{cond2}) essentially shifts the weight on control to
the second term of the target function (\ref{sample_target}) so that the
function asymptotically takes the form%
\begin{equation}
P_{k}\simeq \frac{\lambda _{0}}{2h}(\delta u_{k})^{2}\text{, \ \ as \ \ }%
h\rightarrow 0  \label{target4}
\end{equation}

Such a target function leads to the solution $\delta u_{k}=0$, which
effectively eliminates the control equation. In other words, the iterative
algorithm seems to be essentially discrete. As a result, the control input $%
u_{k}$, generated by the first equation in (\ref{sample_iter1}), depends
upon sampling time interval $h$. Let us illustrate this observation by
implementing the iterations (\ref{sample_iter1}) under the fixed set of
parameters, $\varepsilon =0.1$, $\Omega =1.0$, $\alpha =1.5$, $\zeta =0.025$%
, and initial conditions, $u_{0}=0.025$, $x_{1,0}=1.0$, $x_{2,0}=0.1$, $%
v_{1,0}=v_{2,0}=0$. The values to vary are two different sampling time
intervals, $h=0.01$ and $h=0.001$, and three different values of the damping
constant, $\lambda =0.1$, $\lambda =1.0$, and $\lambda =10.0$. For
comparison reason, Fig. 1 shows time histories of the system coordinates
under the fixed control variable $u=\zeta $. This (no control) case
corresponds to free vibrations of the model (\ref{nb3}) whose dynamics
represent a typical beat-wise decaying energy exchange between the two
oscillators. As mentioned at the beginning of this section, the beats are
due to the 1:1 resonance in the generating system ($\varepsilon =0$, $%
u=\zeta =0$); more details on non-linear features of this phenomenon, the
related analytical tools, and literature overview can be found in \cite%
{Pilipchuk:2010} and \cite{ManevitchGendelman:2011}. In particular, the
standard averaging method was applied to the no damping case of system (\ref%
{nb3}) in \cite{Pilipchuk:2010}.

Now the problem is to suppress the beat phenomenon by preventing the energy
flow from the first oscillator into the second oscillator. As follows from
Figs. 2 through 5, such a goal can be achieved by varying the damping ratio
of the second oscillator, $\{u_{k}\}$, during the vibration process
according to the algorithm\footnote{%
Note that, although the algorithm is designed to suppress accelerations of
the second oscillator, acceleration and energy levels of vibrating systems
are related.} (\ref{sample_iter1}). First, the diagrams in Figs. 2 and 3
confirm that the sampling time interval $h$ represents an essential
parameter of the entire control loop. In particular, decreasing the sampling
interval from $h=0.01$ to $h=0.001$ effectively increases the strength of
the control; compare fragments (b) in Figs. 2 and 3. However, if such
decrease of the sampling time is accompanied by the increase of $\lambda $
according to condition (\ref{cond2}), then the strength of control remains
practically unchanged; compare now fragments (b) in Figs. 2 and 4. As
follows from fragments (a) in Figs. 2 and 4, the above modification of both
parameters, $h$ and $\lambda $, also brings some difference in the system
response during the interval $80<t<150$, but this is rather due to numerical
effect of the time step. Finally, analyzing the diagrams in Figs. 3 and 5,
shows that reducing the parameter $\lambda $ as many as ten times under the
fixed time step $h$ leads to a significant increase of the control input $%
\{u_{k}\}$ with a minor effect on the system response though. Therefore the
parameter $\lambda $ can be used for the purpose of satisfying some
constraint conditions on the control inputs $\{u_{k}\}$ in case such
conditions are due to physical limits of the corresponding actuators. In
addition, let us show that parameter $\lambda $ may affect the convergence
of algorithm (\ref{sample_iter1}) based on the following convergence
criterion \cite{Ostrowski:1960}:

\textit{For a fixed point }$z_{\ast }$\textit{\ to be a point of attraction
of the algorithm }$z_{k+1}=G(z_{k})$\textit{\ a sufficient condition is that
the Jacobian matrix of }$G$\textit{\ at the point }$z_{\ast }$\textit{\ has
all its eigenvalues numerically less than 1, and a necessary condition is
that they are numerically at most 1. The geometric rate of convergence is
the numerically largest eigenvalue of this Jacobian.}

Applying this criterion to the algorithm (\ref{sample_iter1}) at zero point,
gives that one of the eigenvalues is always zero, $q_{0}=0$, whereas another
four eigenvalues, $q_{i}$ ($i=1,...,4$) are proportional to the time step, $%
q_{i}=hp_{i}$, where the coefficients $p_{i}$ are given by the roots of
algebraic equation%
\begin{equation}
p^{4}+2\zeta \Omega p^{3}+2\Omega ^{2}p^{2}+2\zeta \Omega
^{3}p+(1-\varepsilon ^{2})\Omega ^{4}=0  \label{char equation}
\end{equation}

As follows from (\ref{char equation}), the damping coefficient $\lambda $
has no influence on the convergence condition near the equilibrium point,
and the convergence can always be achieved under a small enough time step $h$%
. Nevertheless, the damping coefficient may appear to affect the convergence
away from the equilibrium point. In this case, analytical estimates for
eigen values of the Jacobian become technically complicated unless $%
\varepsilon =0$, when four of the five eigenvalues vanish as $h\rightarrow 0$%
, except one eigenvalue, which is estimated by 
\begin{equation}
q=-\left( 1+\frac{\lambda }{4\Omega ^{2}v_{2}^{2}}\right) ^{-1}  \label{q}
\end{equation}

This root gives $q\rightarrow q_{0}=0$ as $v_{2}\rightarrow 0$. However,
when $v_{2}\neq 0$, equation (\ref{q}) gives the estimate $0<q\leq 1$ as $%
\infty >\lambda \geq 0$. Therefore, only the necessary convergence condition
is satisfied for $\lambda =0$.

\section{Conclusions}

In this work, a discrete time control algorithm for nonlinear vibrating
systems using the damped least squares is introduced. It is shown that the
corresponding damping constant $\lambda $ and sampling time step $h$ must be
coupled by the condition $\lambda h$ = \textit{constant} in order to
preserve the result of calculation when varying the time step. In
particular, the above condition prohibits a direct transition to the
continuos time limit. This conclusion and other specifics of the algorithm
are illustrated on the nonlinear two-degrees-of-freedom vibrating system in
the neighborhood of 1:1 resonance. It is shown that the dissipation ratio of
one of the two oscillators can be controlled in such way that prevents the
energy exchange (beats) between the oscillators. From practical standpoint,
controlling the dissipation ratio can be implemented by using devices based
on the physical properties of magnetorheological fluids (MRF) \cite%
{Dyke:1996}, \cite{Phule:2001}. In particular, different MRF dampers are
suggested to use for semi-active ride controls of ground vehicles and
seismic response reduction.

\newpage

\includegraphics[scale=0.75]{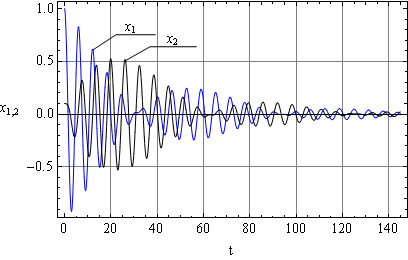}

Figure 1: No control beat dynamics with the decaying energy exchange between
two Duffing's oscillators; $u=\zeta =0.025$.

\includegraphics[scale=0.75]{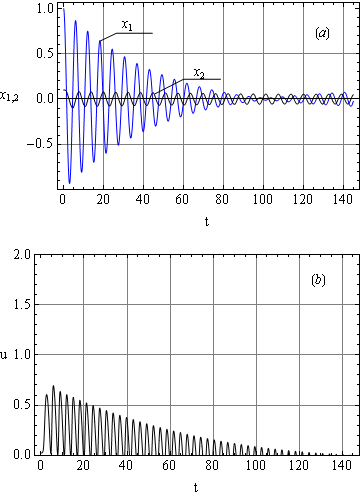}

Figure 2: Beat suppression under the time increment $h=0.01$ and weight
parameter $\lambda =1.0$: (a) the system response, (b) control input - the
damping ratio of second oscillator.

\includegraphics[scale=0.75]{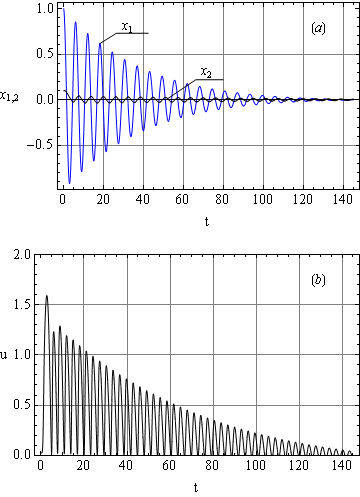}

Figure 3: Beat suppression under the reduced time increment $h=0.001$ and
the same weight parameter $\lambda =1.0$: (a) the system response, (b)
control input - the damping ratio of second oscillator.

\includegraphics[scale=0.75]{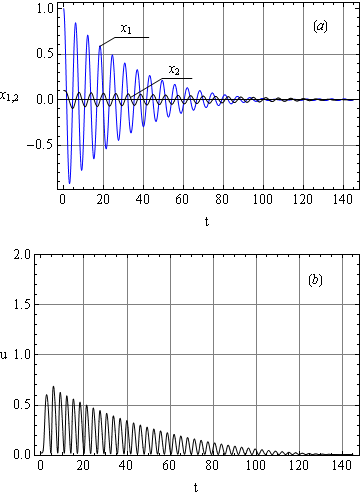}

Figure 4: Beat suppression under the reduced time increment $h=0.001$ but
increased weight parameter $\lambda =10.0$: (a) the system response, (b)
control input - the damping ratio of second oscillator.

\includegraphics[scale=0.75]{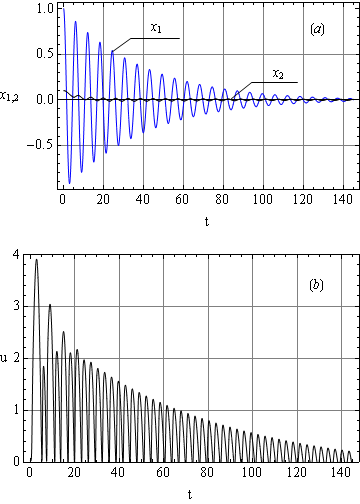}

Figure 5: Beat suppression under the reduced time increment $h=0.001$ and
vanishing weight parameter $\lambda =0.1$: (a) the system response, (b)
control input - the damping ratio of second oscillator. 

\end{document}